\input amstex
\documentstyle{amsppt}
\magnification\magstep1
\nologo
\NoBlackBoxes
\NoRunningHeads
\pagewidth{28pc}
\pageheight{43pc}
\topskip=\normalbaselineskip
\def\Int{\mathop{\roman{Int}}\nolimits}
\def\Cl{\mathop{\roman{Cl}}\nolimits}
\def\oo{\varnothing}
\let\ge\geqslant
\let\le\leqslant
\def\C{{\Bbb C}}
\def\R{{\Bbb R}}
\def\Z{{\Bbb Z}}

\def\Cp#1{\C\roman P^{#1}}
\def\CP#1{\C\roman P^{#1}}
\def\Rp#1{\R\roman P^{#1}}
\def\barRP#1{\overline{\R\roman P}^{#1}}
\def\barCP#1{\overline{\C\roman P}^{#1}}
\def\e{\varepsilon}

\def\dsum{\bot\!\!\!\bot}
\def\ra{\rangle}
\def\la{\langle}
\def\a{\alpha}
\def\b{\beta}
\let\tm\proclaim
\let\endtm\endproclaim

\let\rk=\remark
\let\endrk=\endremark

\def\conj{\mathop{\roman{conj}}\nolimits}
\def\c#1{\conj^{#1}}
\def\W#1{W^{#1}}
\def\X#1{X^{#1}_{\R}}
\def\Y#1{Y^{#1}}
\def\p#1{p^{#1}}
\def\A#1{\goth A^{#1}}
\def\red{\succ}
\def\wed{\vartriangleright}
\def\shh#1#2{\langle{#1}\dsum 1\langle{#2}\rangle\rangle}
\def\shhh#1#2#3{\shh{#1}{#2}_{#3}}
\def\ss#1#2{\langle{#1}\rangle_{#2}}
\def\sss{\langle1\langle1\langle 1\rangle\rangle\rangle}
\def\shhhh#1#2#3#4{\shhh{#1}{#2}{#3}^{#4}}

\def\aa#1{\goth A_{#1}}
\def\xrr#1{X_{\R, #1}}

\def\bn#1{b_2^{#1}(Y)}
\def\bp#1{b_2^+(#1)}
\def\bm#1{b_2^-({#1})}
\def\cil{ Z}
\topmatter
\title
Rokhlin Conjecture and
Topology of Quotients
of Complex Surfaces by Complex Conjugation
\endtitle
\author
S\. M\. Finashin
\endauthor
\address
Middle East Technical University,
Ankara, 06531, Turkey
\newline\indent
St\.-Petersburg Electrotechnical University,
199376, Russia
\endaddress
\email serge\,\@\,rorqual.cc.metu.edu.tr
\endemail
\abstract
Quotients $Y=X/\conj$ of complex surfaces by anti-holomorphic
involutions $\conj\: X\to X$ tend
to be completely decomposable when they are simply connected,
i.e., split into connected sums, $n \Cp2\#m\barCP2$, if $w_2(Y)\ne0$,
or into $n(S^2\times S^2)$ if  $w_2(Y)=0$.
If $X$ is a double branched covering over $\Cp2$, this phenomenon
is related to unknottedness of Arnold surfaces in $S^4=\Cp2/\conj$,
which was conjectured by V.Rokhlin. The paper contains proof of
Rokhlin Conjecture and of decomposability of quotients for
plenty of double planes and in certain other cases.

This results give, in particular, an elementary proof of
Donaldson's result on decomposability of $Y$ for $K3$ surfaces.
\endabstract
\endtopmatter
\document

\heading
\S 1. Introduction
\endheading

By a Real variety (e.g., Real surface, Real curve)
we mean a pair $(X,\conj)$, where
$X$ is a complex variety and $\conj\:X\to X$ 
an anti-holomorphic involution called
{\it the real structure} or {\it the complex conjugation}.
Given an algebraic variety over $\R$
we consider its complexification with the natural complex
conjugation (the Galois transformation) as the corresponding
Real variety.

 This paper is devoted to studying the topology of
quotients $Y=X/\!\conj$ for nonsingular Real surfaces $(X,\conj)$.
 The quotient $Y$ inherits from $X$ an orientation and a
smooth structure, which
makes the projection $q\: X\to Y$ an orientation preserving and smooth
$2$-fold covering branched
along the real part $X_{\R}$ of $X$,  the fixed point set of \ $\conj$,
which is identified, here and in what follows,
with its image $q(X_\R)$.

 It turns out that $Y$ tends to split
into a connected sum of elementary pieces,
$\CP2$, $\barCP2$ and $S^2\times S^2$,
when it is simply connected;
in particular, when $X$ is simply connected and $X_{\R}\ne \oo$.
We call this property of complete decomposability for quotients
CDQ-property and call
$(X,\conj)$ a CDQ-surface 
if it is satisfied.
CDQ-property provides a link between
the differential topology
  of complex surfaces and  differential topology of
knotted surfaces in the decomposable $4$-manifolds,
$X_{\R}\subset Y$, which was used in \cite{FKV}
in construction of exotic knottings.
A plenty of examples may suggest that
CDQ-property is universal for Real surfaces, i.e. holds 
whenever $Y$ is simply connected;
we will refer to it as to CDQ-conjecture.

The first and rather famous example of CDQ-surface
 is $\Cp2$ with
 $\Cp2/\!\conj\cong S^4$ ($\cong$ is read ``diffeomorphic''),
which is known as
 Massey--Kuiper theorem.
This fact has a long history as a folklore
(see \cite{A1, A2}),
but, its proof, perhaps, was not published before \cite{Ma, K}.
The next examples 
are rational surfaces,  $X\cong \CP2\,\#\, n\barCP2$,
obtained by  blow-ups at real points of $\Cp2$.
The Massey--Kuiper theorem prompts that
the quotient will be the same,
$Y\cong S^4\,\#\, n S^4\cong S^4$.
CDQ-property for
one more class of Real rational surfaces was set up in \cite{Le}.

Logarithmic transforms of an elliptic surface
can also be made in the the real category (real fibers on a Real
surface) and does not change the quotient as well \cite{FKV}.
For example, by logarithmic transforms on
$E(1)\cong \CP2\, \#\, 9\barCP2$ 
we can get CDQ-Dolgachev surfaces, $D_{p,q}$, and, moreover,
 Real elliptic surfaces with arbitrary number of
multiple fibers and arbitrary multiplicities
with the quotient $Y$ being still diffeomorphic to $S^4$.
Further, S.Donaldson noticed in \cite{D} that CDQ-property
holds for K3 surfaces and
S.~Akbulut  \cite{Ak}
gave examples of CDQ-surfaces of general type by
 proving CDQ-property 
for a series of Real double planes.

We prove CDQ-property
for several new families of Real surfaces,
which include a plenty of double planes
and, more generally, 
certain doubles of CDQ-surfaces.
Some related problems
of topology of algebraic curves
and the topology of $X_{\R}$ in $Y$
are also studied.

In $\S 2$ we 
introduce Arnold surfaces for Real algebraic curves and
discuss the Rokhlin Conjecture, which is in fact
a refined version of the CDQ-conjecture
 in the case of double planes.
In $\S3$ we prove 
the Rokhlin Conjecture for curves which can be
 obtained by  perturbation
of  lines in generic position.
In $\S4$  we discuss some phenomena related to deformations of
Real curves. The effect of such deformations for $Y$ is well-known
\cite{Le, W}. We give a more subtle description, which takes
into account the topology of $X_\R$ in $Y$ and the topology of
  Arnold surfaces in $S^4$.
We show
 that vanishing of a torus component $T\subset X_\R$
leads to a logarithmic transform of $Y$ along $T$.
In $\S5$ we prove a generalized version of the Rokhlin Conjecture
for a certain class of ``double'' curves in CDQ-surfaces.
As a corollary we get unknottedness for the images $q(A)\subset Y$
 of imaginary curves $A\subset X$.
In $\S6$
a special case of Generalized Rokhlin Conjecture and
Decomposability
Problem is considered
for fibered Real surfaces.
In the case of elliptic surfaces, $P=E(n)_{m_1,\dots,m_k}$,
our results imply
that the exotic phenomena which were found in \cite{FKV} for
 the  knottings $P_{\R}\subset Q=(P/\!\conj)$,
cannot straightforwardly
produce  exotic differential structures in the quotients
of the double coverings over $P$.

In $\S7$ we discuss some applications
for Real plane curves of degree $\le6$ and
give in particular an elementary proof of
CDQ-property for Real K3 surfaces,
alternative to the Donaldson's.

It may be worth to mention some other recent related results.
The author proved CDQ-property
for {\it all\, } Real rational and Enriques surfaces \cite{F3}
and for {\it certain\/} complete intersections of {\it arbitrary\/}
multi-degree, namely, for the ones which can be
constructed by method of a small perturbation \cite{F2}.
Note also that vanishing of Donaldson  (Seiberg--Witten) invariants for $Y$,
is a weaker version of CDQ-property,
which is often easier to prove.
For example, vanishing in the case when
$X_\R$ contains an orientable component of genus $\ge 2$
is a trivial corollary of the adjunction formula
  (see \S4).
Another example: if $X_\R=\oo$ then
$Y$ is undecomposable, however, its
   Seiberg--Witten invariants vanish by the recent result of S.Wang.

\subheading{Acknowledgments}
I would like to thank S.~Akbulut for several helpful references
and remarks.

\heading
\S 2. Rokhlin Conjecture
\endheading
Given  a nonsingular curve $ A\subset \Cp2 $,
the zero set of a degree 2k real form $f$,
consider the surfaces
$\W{\pm}=\W{\pm}(A)=\{x\in \Rp2 \ | \ \pm f(x) \ge 0 \}$.
We have obviously $\W+ \cup \W-= \Rp2,\ \W+ \cap  \W-= A_{\R} $.
One of these surfaces is, obviously,
orientable, the other is not.
Being free to change the sign of $f$ we can assume that the orientable
surface is $\W+$.
By taking unions,
$\A{\pm}=\A{\pm}(A)=(A/\!\conj)\cup\W{\pm}$,
we get 2 closed surfaces in
$ S^4=\Cp2/\!\conj$
known as {\it Arnold surfaces}
in Topology of Real Algebraic Curves.
They are not smooth along $A_{\R}$ but can be smoothed,
since $A$ intersects $\Rp2$ normally.
In what follows
we will leave the notation $\A{\pm}$ for
the smoothed Arnold surfaces.

\proclaim{2.1. Rokhlin Conjecture }
  Arnold surfaces are standard for nonsingular Real curves $A$
 of even degree with $A_{\R}\ne\varnothing$.
\endproclaim

By a standard surface in $S^4$ we mean a connected sum of
several copies of standard tori, or  standard  $\Rp2$.
A standard torus, which can be characterized as a boundary
of a solid torus embedded in $S^4$, is unique up to isotopy.
 $\Rp2$ can be embedded in $S^2$ in two (up to isotopy) standard ways
distinguished by the normal numbers, which can be  $\pm2$.
We will denote a standard torus by $(S^4, T^2)$ and standard $\Rp2$
with the normal numbers $-2, +2$
by $(S^4,\Rp2)$
and  $(S^4,\barRP2)$ respectively.

Note that surfaces which bound
 solid handlebodies embedded in $S^4$
are standard.
(By a solid handlebody we mean a 3-manifold homeomorphic
to $F\times[0,1]$, where $F$ is a connected compact
surface with $\partial F\ne\oo$.)

Note also that $(\Cp2,\Rp2)/\!\conj=(S^4,\Rp2)$.

\rk{Remarks}
\roster
\item
When $A_{\R}=\oo$ the Rokhlin Conjecture still may concern
$\A+=A/\!\conj$. On the other hand, $\A-$
has two  components, $\Rp2$ and   $\A+$, which are linked, since
the double covering over $S^4$ branched along  $\A-$ has
fundamental group $\Z/2$ 
 (it would be $\Z$
for an unlink).
\item
There is a weaker version of the Rokhlin Conjecture studied in~\cite{F1}:
\newline
{\it
Complements of the Arnold surfaces
for nonsingular Real curves $A$  of even degree, with
$A_{\R}\ne\varnothing$,
have abelian fundamental groups.
}
\endroster
\endrk

We define in addition Arnold surfaces in $\Cp2$
when $A-A_{\R}$
splits in two connected components, which happens 
 if the fundamental class
$[A_{\R}]\in H_1(A;\Z/2)$ vanishes. 
If this is the case  $A$ is said to be of
 type $1$, otherwise it is
of  type $2$.
If $A$ is of type $1$ then we denote by $A_k, k=1,2$
the components of $A-A_{\R}$, put
$\A{\pm}_k=\W{\pm}\cup A_k$ and call $\A{\pm}_k$
the Arnold surfaces in $\Cp2$.
 The Rokhlin conjecture may treat them
if we define a standard surface
in $\Cp2$ as a surface which is
isotopic to a connected sum
of a nonsingular algebraic curve in $\Cp2$ with
a standard surface in $S^4$.

The Rokhlin conjecture for $\A\pm(A)$ is stronger then CDQ-conjecture for
 double planes $X$ branched along curves $A\subset\Cp2$,
as it follows from the following construction
well known after
\cite{A1}.

Define $X$ by the equation
$f(z_0\!:\!z_1\!:\!z_2)=w^2$
in the weighted quasi-projective 3-dimensional space
with coordinates $z_0,z_1,z_2$ of weight $1$
and $w$ of weight $k$.
It has 2 real structures defined by
the covering conjugations:
$\conj^{\pm}(z_0:z_1:z_2:w)=
(\overline{z_0}:\overline{z_1}:\overline{z_2}:\pm\overline{w}).$

The real part $\X{\pm}$ of $(X, \conj^{\pm})$
is mapped onto $\W{\pm}$
by the projection $p\:X\to \Cp2$
as a double covering
with sheets glued together along $\partial\W{\pm}=A_{\R}$.
For even $k$ this covering
is trivial, hence, $\X{\pm}$
is the usual double of $\W{\pm}$.
For odd $k$ it is the orientation
covering, so,
  $\X{\pm}$ is orientable.

Further,
$\conj^+$ and $\conj^-$ commute,
give in product the covering transformation
and thus define an action of $\Z/2\oplus\Z/2$ on $X$.
Factorization by $\Z/2$ in  different order
yields the following commuting diagram:
$$
\CD
X     @>q^{\pm}>>   \Y{\pm}=X/\!\conj^{\pm}\\
@VpVV         @VV\p{\pm}V\\
\Cp2   @>q>>   S^4=\Cp2/\!\conj  ,
\endCD
$$
where $ q, q^{\pm}$ are the quotient maps and
$\p{\pm}$ is the double covering branched along
$\A{\mp}$ (cf.~\cite{A1}).

This description implies

\proclaim {2.2. Theorem}
If $\A{\pm}$ for a curve $A$ is standard, then
$(X,\conj^\mp)$ is CDQ-surface.
\endproclaim

\demo{Proof}
Use
 the above diagram  and note that
the double branched coverings over
$(S^4, T^2)$, $(S^4, \Rp2)$ and $(S^4,\barRP2)$
are $S^2\!\times\! S^2$, $\Cp2$ and $\barCP 2$ respectively,
and that
the double covering over
$(Y_1\#Y_2,F_1\#F_2)$
 is the connected sum of the double coverings over
$(Y_1,F_1)$ and $(Y_2,F_2)$.
\qed
\enddemo

\heading
\S 3.  Rokhlin Conjecture for $L$-curves
\endheading

Recall
that the parameter space $C_m$ of Real plane algebraic curves
of degree $m$ is a real projective space
of dimension $\binom{m+2}2 -1$ and
singular curves constitute the discriminant hypersurface,
$\Delta_m\subset C_m$.
A {\it deformation} $A_t$ of $A_\alpha\in C_m$
is a path $[\alpha,\beta]\to C_m$, $t\mapsto A_t$.
A {\it perturbation} of $A_\alpha\in \Delta_m$ is a deformation $A_t$
such that $A_t\in~C_m-\Delta_m$, $\forall~t\in(\alpha,\beta]$.

We call $A_\beta$ an {\it $L$-curve} if it can be obtained
by a perturbation of a curve $A_\alpha$ which splits into a union
of $m$ Real lines in general position (three lines should not
have a common point).

\tm{3.1. Theorem\!}
The Rokhlin Conjecture holds for $\A{\pm}(A)\!\subset\! S^4$
for $L$-curves~$A$.
\endtm

\tm{3.2. Corollary }
Double planes $X\to \Cp2$ branched along $L$-curves
are CDQ-surfaces.
\endtm

\demo\nofrills{ Proof }\ of Theorem 3.1
 is based on the following criterion~:

\tm{3.3. Theorem \cite{Li}}
Let $F\subset D^4\subset S^4$ be a closed surface
which lies in $S^3=\partial D^4$, except  several
$2$-discs  standardly embedded inside 4-disc $D^4$.
Then
$F$ is a standard surface in $S^4$.
\endtm

Let us describe first a suitable handle decompositions
 of $\Cp2$ and $\Cp2/\!\conj$.
Choose a point $b_0\in \Rp2$ with a small regular $\conj$-symmetric
neighborhood $B_0$ around it. Consider the central projection
$r\: (\Cp2-B_0)\to C$, from $b_0$ to some Real line $C$, $ b_0\notin C$.
Let $C_1, C_2,$ be the closures of the connected components of $C-C_{\R}$
and $B_i=r^{-1}(C_i),\ i=1,2$.
Then
$\widetilde{H}=B_1\cap  B_2$ is a solid torus, a disc bundle over $C_{\R}$,
which contains a M\" obius band $M=\Rp2-B_0$.
4-discs $B_1$, $B_2$, $B_0$ can be considered as
 0-, 2-  and 4-handles giving a decomposition
of $\Cp2$, like in \cite{Ak}.
The quotient $q(B_0)$ is again a 4-disc.
$B_1$ and $B_2$ are permuted by $\conj$, hence,
$B=q(B_1\cup B_2)= B_1/\sim$,
where $\sim$ identifies conjugated points in ${\widetilde H}$,
is a 4-disc as well. $H=\widetilde H/\!\conj$ can be thought of
as the trace of an isotopy of $M$, fixed on $\partial M$,
which pushes $M$ outside of the interior,
$\Int(B)$, to $\partial B$.

Consider now an $L$-curve $A$ obtained from
$A_0=L_1\cup\dots\cup L_m$ by a perturbation.
It is well known that a
 perturbation of a curve does not change its isotopy type
outside a neighborhood of singularities. 
So,
$A$ can be assumed
to coincide with $A_0$
in the complement of small regular $\conj$-symmetric open discs
$\widetilde P_{ij}$
around points $ p_{ij}=L_i\cap  L_j, \ 1\le i< j\le m$.
Put $P_{ij}=q(\widetilde P_{ij})$ and
$ P=\bigcup_{1\le i<j\le m}P_{ij}$.
We can assume also that
$b_0$ is chosen
inside $\W{\mp}$ and that $B_0$ does not intersect  $A_0$ and $P_{ij}$.

The idea behind the proof is
to apply Theorem 3.3 to $\A\pm$ using
the 4-disc obtained from $B$ by cutting it along $H$.
The result can be
 naturally identified with $B_1$ (or with $B_2$).
It contains $\W{\pm}$ in its boundary and $(A/\!\conj)-P$ consists
of $m$ standard discs $(L_i/\!\conj)-P$.
More formally, we get a 4-disc $B'$ after deleting from $B$ a small regular
neighborhood $N$ of $H-M$.
Then we push by an isotopy the part of $\A{\pm}$
contained inside $N$
to $\partial B'$.

It suffices to construct such an isotopy inside $P$,
since $(\A{\pm}-P)\cap \Int(N)=\varnothing$.
So, we only need to study the local question,
specifically,
possible  positions of $\A{\pm}$ with respect to $H$
inside  $P_{ij}$.
It is determined, obviously,
 by the position of the lines $L_i$, $L_j$, with respect
to the basepoint $b_0$ and with respect to
$W^{\pm}$. The 4 possible cases are shown
on Figure 1.
\midinsert
\topcaption{Figure 1}
Positions of $L_i, L_j$
with respect to $ \W{\pm}$ (the shaded region)
and to $b_0$
\endcaption
\vspace{40mm} 
\endinsert
Let
$W_{ij}=\W{\pm}\cap P_{ij}$, $A_{ij}=(A/\!\conj)\cap P_{ij}$, 
$H_{ij}=H\cap
P_{ij}$,
and $\A{}_{ij}=\A{\pm}\cap P_{ij}$.
\ $\partial A_{ij}$ splits 
into 2 pairs of arcs $l_1=A_{ij}\cap\partial P_{ij}$
and $l_2=A_{ij}\cap W_{ij}$.
By standard arguments of \cite{M}, $A_{ij}$
can be pushed out to $\partial P_{ij}$ by an isotopy which is
transversal to $H$,  keeping fixed on $l_1$
 and moves $l_2$ in $\Rp2$ so that $A_{ij}$ is preserved
not to intersect normally with the latter.
It gives an isotopy of $\A{\pm}$ which is fixed on $(A/\!\conj)-P$,
preserves $\W{\pm}-P$ invariant and  pushes $A_{ij}$ to $\partial P_{ij}$.
In the cases (a)  and (b)
this isotopy pushes $\A\pm$ outside of the interior of $P_{ij}$.
In the case (d),\  the part $W_{ij}$ is still inside of $P_{ij}$
but it can be
pushed out to $\partial P_{ij}$ along $H_{ij}$
by an isotopy
without obstructions (see Figure 2).
\midinsert
\topcaption{Figure 2}
 Mutual position of $H$ (the disc) and $A_{ij}$ (the band)
after pushing  $A_{ij}$ off $\Int(P_{ij})$.
$\A{}_{ij}$ is shaded 
\endcaption
\vspace{30mm}
\endinsert

In the case (c)  we cannot do the same, because
$A_{ij}$ is intersecting  $H$ not only at
the part of its boundary $l_2$ but also
along the ``middle line'' $l_3$
(see Figure 3).
\midinsert
\topcaption{Figure 3}
$l_3\subset A_{ij}\cap H$ 
formed by intersection of $A$ with the pencil of Real lines
passing through $b_0$.
The isotopy leaves
 $W_{ij}'$ (which is shaded on the right)
 inside $P_{ij}$
\endcaption
\vspace{35mm}
\endinsert
So,  after being pushed off $\Int(P_{ij})$  
along $H$, $A_{ij}$ intersects $W_{ij}$
as it is shown on Figure 2.
 But, we can avoid this  and
 push to $\partial P_{ij}$  only a
 neighborhood of $l_1$
 leaving the rest, the neighborhood $W_{ij}'$ of $l_3$, inside $P_{ij}$.

Finally,  we add to $B'$ the discs $P_{ij}$ for which $\A{}_{ij}$
lie in the positions (a) and (c) and 
excise from $B'$ the other ones.
We get, clearly, a new  4-disc, $B''$ and can apply now Theorem 3.3
to $\A\pm$.
To see how  $\A\pm$ lies in $B''$ after this isotopy
note that in the cases (a) and (c)
$\partial B''\cap\partial P_{ij}=\Cl(N\cap \partial P_{ij})$
is a proper neighborhood of the 2-disc $H\cap\partial P_{ij}$.
It is clear from Figure 2 that we can
vary  size of this neighborhood by varying size of $N$
so that it will contain $A_{ij}$ in the cases (a),
and both $A_{ij}$  and
$W_{ij}-W_{ij}'$  in the case (c).
In the cases (b) and (d)
$\partial B''\cap\partial P_{ij}$
is the complement of the neighborhood
$N\cap \partial P_{ij}$ in $\partial P_{ij}$.
It will contain $\A{}_{ij}$ in the case (b)
if the neighborhood $N$ is small enough. It will also contain
$\A{}_{ij}$ in the case (d) 
 after we push
$W_{ij}$ to $\partial N\cap \partial P_{ij}$ by an isotopy.

After the above isotopy
$\A{\pm}$ lies
in the boundary of $B''$,  except  $m$ discs
$(A/\!\conj)\cap \Int(B'')$
left from $L_i$
and discs $W_{ij}'$,
which lie all inside $B''$  and are obviously unknotted.
\qed
\enddemo

\rk{Remarks}
\roster
\item
$L$-curves form a pretty large collection of Real curves.
They include, in particular, ``maximal nest'' curves,
that is degree $2k$ curves
which real parts consist of $k$ ovals linearly ordered by inclusion
(see the example on Figure 8).
Thus Theorem 3.1 generalizes results of \cite{Ak},
where maximal nests were under consideration.
\item
One can  easily notice the upper bound, $\frac13m(m-1)$,
for the number of real components of $L$-curves of degree $m\ge3$.
For $m\ge 5$
it is less
then the Harnak bound, $\frac12(m-1)(m-2)+1$, effective
for arbitrary Real curves of degree $m$.
This shows one of the restrictions to real schemes of $L$-curves.
More examples  and discussions of $L$-curves can be found in \cite{F1},
\cite{FKV} (see also \S7 below).
\endroster
\endrk

\heading
\S4.
Deformations of real structures
\endheading
In this section we describe how the
topology of $\A{\pm}$, $Y$  and $X_{\R}$
is changing along with the
 deformation of
the branched locus $A$.

Any pair of nonsingular plane  Real curves,
$A_\alpha, A_\beta\in C_m$,
 can be
connected by a deformation, say $A_t$, $t\in[\alpha, \beta]$.
If $A_t\notin\Delta_m$, $\forall t\in[\alpha, \beta]$,
the deformation is called
{\it rigid isotopy}.
In this case the $\conj$-equivariant
topological type of $(\Cp2, A_t)$ and, hence,
the topological type
of $\A{\pm}$, $Y$, $X_{\R}$
remains unchanged.
If $A_t$ is a generic deformation then it
  crosses $\Delta_m$ transversally at several nonsingular points.
Such a nonsingular point, say, $A_0\in\Delta_m$,
represents a curve
with one real node,
$P\in A_0$, which can be
either isolated in $\Rp2$
or cross-like.
The topological effect
of transversal crossing of $\Delta_m$
for $\A{\pm}$, $Y$, $X_{\R}$, is local and
 depends, obviously, only on the type of the node at $P$
and on the local position of $P$ with respect to $W^{\pm}$.
Figure 4 shows $6$ possible topologically different elementary
 modifications of $W^\pm$, which will be called moves of $W^\pm$
and denoted by $M_i$ or $M_i^{-1}$, $i=0,1,2$.
\midinsert
\topcaption{Figure 4}
Elementary moves.  $W$ is the shaded region
\endcaption
\vspace{45mm}
\endinsert

Consider a variation $A_t, t\in [-1,1]$ which crosses $\Delta_m$
transversally at $t=0$.  
Denote by $X_t$ the double plane branched along $A_t$ and by
$W_t$ one of the
domains $\W\pm(A_t)$
 continuously following the curves $A_t$
(we omit the superscript
sign over $W_t$,  since the move $M_1$ may change it).
Denote by
$\aa t$, $\conj_t$, $ X_{\R,t}$, $Y_t$ the
objects which follow continuously
  $W_t$,
so that if $W_t=\W\pm(A_t)$, then
$\A{}_t=\A\pm(A_t)$, but $\conj_t=\conj^\mp$,  $ X_{\R,t}= X_{\R,t}^\mp$
 and $Y_t=Y^\mp$. 
Assume that $M$ is the move which describes
 the modification of $W_t$ as $t$ increases.

\tm{4.1. Theorem} If $M=M_1$ or $M_0^{-1}$ then
\roster
\item
$\aa 0$ is smoothable at $P$,
hence, $Y_0$ is smoothable as well;
\item
$(S^4,\aa{0})\cong(S^4,\aa{-1})$, hence, $Y_0\cong Y_{-1}$;
\item
$(S^4,\aa{1})\cong(S^4,\aa 0)\#(S^4,\barRP2)$, hence,
$Y_{1}\cong Y_0\#\barCP2$;
\item
 after smoothing of $Y_0$, the real part,
$\xrr 0\subset Y_0$,
is locally diffeomorphic at $P$ to a node,
i.e., there exist an orientation preserving diffeomorphism $U\cong\C^2$
of a sufficiently small regular neighborhood $U\subset Y_t$ of $P$
 which sends $X_{\R,0}\cap U$ to the union of coordinate
lines;
\item
$(Y_{-1},\xrr{-1})$ is obtained from $(Y_{0},\xrr{0})$ 
by an ``algebraic'' perturbation 
of the above singularity (i.e., the local effect of this deformation
in the  chart $U\cong\C^2$ is the perturbation of the node);
\item
$(Y_1,\xrr1)$ is obtained from $(Y_0,\xrr0)$ by blowing-up at $P$.
\endroster
\endtm

\tm{4.2. Theorem}
Suppose that $M=M_2$. Then
\roster
\item
 $\aa 0$  is locally
diffeomorphic to an algebraic curve
near a simple tangent point singularity;
more precisely,
 there exist an orientation preserving diffeomorphism $U\to\C^2$
of a sufficiently small regular neighborhood $U\subset S^4$ of $P$
 which sends $\aa0\cap U$ to the curve $z_2^2-z_1^4=0$;
in particular, this singularity has $(2,4)$ torus link;
$Y_0$ has a singularity at $P$ which link is the lens space
$L(4,1)$;
\item
$(S^4, \aa{-1})$ is obtained from $(S^4, \aa{0})$ by 
removing a neighborhood of $P$ $($a pair of tangent discs$)$
and replacing it by a Seifert surface
 of the $(2,4)$ torus link 
$($see Figure 5$)$;
it implies that
$Y_{-1}$ is obtained from $Y_0$ by exchanging a neighborhood
of $P$ $($the cone over $L(4,1))$ for
the total space of a $D^2$-bundle over $S^2$ with
the normal number $-4$;
\item
$(S^4,\aa 1)$ is obtained from $(S^4,\aa 0)$ by replacing
 a neighborhood of $P$ by
a disc  and a M\" obius band,
as it is shown on Figure 5;
$Y_{1}$ is obtained from $Y_0$ accordingly,
by substituting of the cone over $L(4,1)$
for the total space of a
 $D^2$-bundle over $\Rp2$ with the normal number $-1$;
\item
$\xrr 0$ contains only an isolated point at $P$ in a neighborhood of $P$;
after the deformation this point disappears in
$\xrr 1$ and gives rise to a sphere component in
$\xrr{-1}$; this sphere can be identified with 
the zero  section of the disc bundle described in (2).
\endroster
\endtm
\midinsert
\topcaption{Figure 5}
Bifurcation of $\aa t$ after the move $M_2$.
The disc component is not seen on the right figure, as it
is pushed inside $D^4$
\endcaption
\vspace{25mm}
\endinsert

\rk{$\bold{4.3.}$ Remarks} $(1)$
The modification of $Y$ after
move $M_2$
 can be thought of as
the ``connected sum'' with $\CP2$ along  a tubular neighborhood
of a conic, say, the Klein conic
defined by the equation
 $z_0^2+z_1^2+z_2^2=0$.
Its complement is  a tubular neighborhood of $\Rp2$.
This sort of surgery
is known as rational blow-down of degree 2
\cite{FS2}.
We can also think of the transform $\aa{-1}\to \aa 1$
as of a ``connected sum'' with $\Rp2$ along a neighborhood of
the circle (real conic)
$x_0^2+x_1^2=x_2^2$ and consider it as a {\it real rational blow-down}.

$(2)$ If instead of perturbation of $A_0$ we consider the resolution,
i.e., blow up $\Cp2$ at $P\in A_0$ and take
the double covering  $\widetilde X_0\to\Cp2\,\#\,\barCP2$ branched along
the proper image $\widetilde A_0$ of $A_0$, further, denote by
$\widetilde{\A{}}_0$ 
 the Arnold surface of $\widetilde A_0$
in $S^4=\Cp2\,\#\,\barCP2/\conj$  (cf.~\S5), by
 $\widetilde X_{\R,0}$ 
the real part of $\widetilde X_0$ and put
$\widetilde Y_0=\widetilde X_0/\conj$,
then
we have $\widetilde{\A{}}_{0}\cong\A{}_{0}$,
$\widetilde Y_0\cong Y_{0}$,
$\widetilde X_{\R,0}\cong X_{\R,-1}$
if the node at $P$ is hyperbolic
(i.e., $M=M_1$ or $M_0^{-1}$),
and
$\widetilde{\A{}}_{0}\cong\A{}_{1}$,
$\widetilde Y_0\cong Y_{1}$, $\widetilde X_{\R,0}\cong X_{\R,1}$
for an elliptic node ($M=M_2$).
\endrk

\tm{4.4. Corollary}
Moves $M_0^{-1}$  and $M_1$ have both the effect of
 a real blow-up on the Arnold surface $($i.e.,~$\#\,\barRP2)$,
 a usual blow-up on $Y$, and a Morse modification
of index 2 on $X_{\R}$.
$M_2$ makes a real rational blow-down on the Arnold surface,
 a rational blow-down of degree 2 on $Y$, and a Morse
modification of index 3 on  $X_{\R}$.
The moves opposite to the above ones have the obvious opposite
effect.
This shows, in particular, how
the modifications on $Y$ is determined by
the Morse modifications of $X_\R$.
\endtm

Given a pair of nonsingular curves, $A_0$, $A_1$, of the same even degree
with 
$W_i=W^\pm(A_i)$,  $\A{}_i=\A\pm(A_i)$, $i=0,1$,
let us write $W_0\red W_1$ if $W_1$ can be obtained from $W_0$
by a deformation connecting $A_0$ and $A_1$ which
involves only the moves $M_0^{-1}$ and $M_1$. 
Clearly, $W_0\red W_1$ implies $\chi(W_0)>\chi(W_1)$.

\tm{4.5. Corollary}
If $\A{}_0$ is standard and $W_0\red W_1$ then $\A{}_1$ is standard
as well.
\endtm

\demo{Proof of Theorems $4.1, 4.2$}
Consider first the case of move  $M_0$.
In a suitable coordinate system near $P$, the projection
$q\:\Cp2\to \Cp2/\!\conj$ is modeled by the map $\C^2\to \C^2$, with
$(z_1,z_2)\mapsto (z_1^2,z_2)$, where $\Rp2$ is represented by
the line $z_1=0$, $A_0$ is defined by $z_1^2-z_2^2=0$.
Then locally $\aa 0=q(A_0)$ is defined by $z_1-z_2^2=0$
and, hence, is smooth and tangent to $\Rp2$ at $P$.
This proves 4.1(1) and 4.2(1).

4.1(2) Near $P$ (i.e.,~inside a small regular neighborhood),
for $t>0$, $\aa t$ is a union of $q(A_t)$ and of
a disc component $D_t$ of $W_t$, bounding an oval
$\partial D_t\subset A_{\R,t}$, which can be considered as
a vanishing cycle for the deformation $t\to 0$.
$\aa 0$ is obtained from $\aa t$ by contracting $D_t$.
Hence, $\aa 0\cong \aa 1$.
4.1(3) For $t<0$, $\aa t=q(A_t)$ is, \ near $P$,  a M\" obius band,
which replaces the disc $D_t$.
The direction of its ``twisting'' is, obviously, standard
and can be easily determined on a model, for example,
conic curve,
by taking double covering branched over it.
4.1(4) Follows from the above local algebraic model describing
the mutual position of $\Rp2$ and $ A_0$ near $P$.
4.1(5) Follows from the arguments of (1) and (2).
$\xrr 0$ is obtained from $\xrr t$ by
the contraction
of $D_t$, topologically
equivalent to the contraction of a vanishing cycle
which gives a double point singularity at $P$.

4.1(6) When  $t\to -0$,
the imaginary vanishing cycle in $\aa t$
is pinched along the disc,
which
is covered by the exceptional curve of the blow-up,
$E\subset Y_t\cong Y_0\#\barCP2$.
$\xrr t\cap E$ consists of a pair of points projected to $P$.
This gives the same topological
description as for the resolution of
a double point singularity.

Let us pass  now to Theorem 4.2.
4.2(2)
A Seifert surface on the $(2,4)$-toric link,
which replaces a pair of tangent discs,
can be easily seen in
 the above local description of $\aa 0$.
4.2(3) The components are a disc part of $\Rp2$ and a M\" obius
band part of $q(A_t)$
(see above 4.1(3)), the mutual position of which is, obviously,
the one described on Figure 5.
4.2(4) The sphere component of $\xrr t$, $t<0$, covers $D_t$;
hence, it is the zero  section. The rest is straightforward.

Now consider the case $M_1$.

4.1(1) The link of $A_0\cup\Rp2$ at $P$ and its quotient
in $\Cp2/\!\conj$ are shown on Figure 6. The link of $\aa 0$
turns out to be an unknot.

\midinsert
\topcaption{Figure 6}
Links for a cross-like node
in $\Cp2$ (on the left) and in $S^4$ (on the right)
\endcaption
\vspace{20mm}
\endinsert

4.1(2)-(3) Can be easily seen from the local description of $\aa t$ shown
on Figure 2.
4.1(4) The link of singularity of $\xrr 0$ is obtained from
the right picture on Figure 6 as the double covering along
the link of $\aa 0$. The result is equivalent
to the left picture.
4.1(5)-(6) Analogous to the case of move  $M_0$, by making use of
the local description of $\A\pm$ on
Figure~2.
\qed
\enddemo

Let us write $W_0\wed W_1$ if $W_1$ can be obtained from $W_0$
by making use of any
moves except $M_2^{-1}$.
It may look plausible that the relation $W_0\wed W_1$
is always true whenever $W_0\ne\Rp2$,
however,  it is not known for curves of degree $\ge 8$.

\smallskip
The following theorems yield
vanishing of Donaldson (Seiber--Witten) invariants for $Y$ in
more general setting,
when the above methods do not allow to prove decomposability of $Y$.
Assume that $b_2^+(Y)>1$ and odd,
which is needed for the invariants to be well defined.
Note  that $b_2^+(Y)=p_g(X)$
(see Supplement 7.4).
It follows from \cite{FS1} and \cite{FS2}
that the vanishing property for Donaldson (Seiber--Witten) invariants 
is preserved by usual blow-ups and -downs and
by rational blow-downs. This immidiately prompts the following

\tm
{4.6. Corollary}
If Donaldson  (Seiber--Witten) invariants  vanish for
 $Y_0$ and $W_0\wed W_1$, then they vanish for $Y_1$ as well.
\endtm

We can also formulate this result in more general setting,
since modifications of $Y$ are determined
by the Morse modifications of $X_\R$.

\tm{4.7. Corollary}
Let $(X_t,\conj_t)$, $t\in[0,1]$  be a generic deformation
of Real surfaces, i.e., $X_t$  nonsingular for all
$t$ with finitely many exceptions for which $X_t$ has one
double point singularity.
Assume that  sphere components of $X_{t,\R}$ do not collapse
while passing double points as $t$ increases.
Then vanishing of Donaldson  (Seiber--Witten) invariants 
 for $X_0/\!\conj_0$ implies their vanishing for $X_1/\!\conj_1$.
\endtm

One can also easily deduce from
 \cite{KM} and \cite{FS2} the following

\tm{4.8. Corollary }
\roster
\item
If $X_{\R}$ contains
an orientable component of genus $\ge2$,
then  the Donaldson polynomials of $Y$ vanish.
\endroster

Furthermore, if $X$ is of simple type then
\roster
\item[2]
fundamental classes of
orientable components of $X_{\R}$ of genus $\ge 1$
are orthogonal in $H_2(X)$
to the Kronheimer--Mrowka  (Seiber--Witten)
basic classes $K_s$ of~$X$;
\item
if $\Sigma\subset X_\R$ is a sphere component, then
$|\Sigma\circ K_s|\le 2$ for any basic class~$K_s$.
\endroster
\endtm
\demo{Proof}
By the adjunction formula (see \cite{KM}, {FS2}),
if $X$ is of simple type, then
$\chi(\Sigma)+\Sigma\circ\Sigma+\max_s(K_s\circ\Sigma)\le0$
for a smooth essential connected surface $\Sigma\subset X$
with $\Sigma\circ\Sigma\ge 0$.
If $\Sigma$ is a component of $X_{\R}$, then
$\chi(\Sigma)+\Sigma\circ\Sigma=0$,
since the normal bundle of  $X_{\R}$ is obtained from the tangent
bundle by multiplication by $i$.
Thus  to prove (2), it suffices to recall that
$J=\max_s(K_s\circ\Sigma)\ge0$, and $J=0$
if and only if $\Sigma$ is orthogonal
to all basic classes $K_s$.
(3) follows similarly from the
adjunction formula for immersed spheres
\cite{FS2}, which  works as well
when $\Sigma\circ\Sigma<0$.

To obtain (1) note that
if $\Sigma\subset X_\R$ has genus
$g(\Sigma)\ge2$ then $\Sigma$ is essential
since it has nonzero self-intersection number.
Furthermore, we have
$\chi(\Sigma)+2\Sigma\circ\Sigma=\Sigma\circ\Sigma>0$,
which contradicts
to the adjunction formula applied to $\Sigma$ in $Y$.
On the other hand, this inequality guarantee
that $Y$ is of simple type.
Hence, the basic classes do not exist.
\qed
\enddemo

By \cite{FS2} a logarithmic transform of multiplicity $m$
can be decomposed into a product of $m-1$ usual blow-ups and
a rational blow-down of degree $m$.
For $m=2$ it is  one blow-up
at a nodal point of a singular
elliptic fiber, which yields a nonsingular rational
curve with self-intersection $-4$, then
a rational blow-down of degree 2 replacing  this curve by $\Rp2$.
 Theorems 4.1---4.2 imply that
  the described combination
 is produced on $Y^{\pm}$ by the combination of
moves $M_1$ and $M_2$ on $W^{\mp}$
which eliminate
 an annulus component of
the complementary domain $W^{\pm}$.

\tm
{4.9. Corollary}
The combination of moves $M_1$ and $M_2$ of $\W{\mp}$
which makes vanish
an annulus component of $\W{\pm}$
as it is shown on Figure 7
produces a logarithmic transform  of multiplicity 2
on $Y^{\pm}$
along the torus component of $X_{\R}$ corresponding to the 
annulus.
\qed\endtm
\midinsert
\topcaption{Figure 7}
Logarithmic transform of multiplicity 2
\endcaption
\vspace{30mm}
\endinsert

Since the modifications involved are local, we can reformulate the above
Corollary as follows.

\tm{4.10. Corollary}
Let  $(X_1,\conj_1)$ be obtained from  $(X_0,\conj_0)$
by a deformation   $(X_t,\conj_t)$,
which consists of nonsingular Real surfaces for all $t\in[0,1]$
with 2 exceptions when $X_t$ has a node.
Suppose that this deformation make 
a torus component of $X_{0,\R}$ vanish.
Then $Y_1$ is obtained from $Y_0$ by a logarithmic
transform of multiplicity $2$ along
this component.
\endtm

\heading
\S 5.
Generalized Rokhlin Conjecture
\endheading

In this section we study
a generalization of the Rokhlin Conjecture,
when $\Cp2$ is replaced by
an arbitrary nonsingular Real surface,
$(P,\conj)$.

Call a curve $A\subset P$ even if it is the zero set of
a holomorphic section  $f: P\to L^{ 2}$, where
 $L^{ 2}=L\otimes L$ is
the square of some holomorphic linear bundle
 $p\: L\to P$.
The latter is called {\it a real bundle} if it
 is supplied with an anti-linear involution
$\conj_{L}\: L \to L$, such that $p\circ\conj_L=\conj\circ p$:
$$
\CD
L     @>\conj{L}>>     L\\
@VpVV              @VpVV\\
P      @>\conj>>       P\\
\endCD
$$

Denote by $L_{\R}$ its real part, that is  the fixed point set of $\conj_L$.
The restriction of $p$ gives a real line bundle
$L_{\R}\to P_{\R}$.
$L^{ 2}$ has the real structure
$\conj_{L^{ 2}}=(\conj_L)^2$ induced from $L$,
with the real part $(L^{ 2})_{\R}$
 trivialized by the choice of the natural positive ray
$(L_{\R})^2$.
Hence,
the sign of $f(x)$ is well defined at real points $x\in P_{\R}$
and we can put
$\W{\pm}=\W{\pm}(f)= \{ x\in P_{\R}  : \pm f \ge 0\} $
 as in the case $P=\Cp2$.
Again, we have
$\W+ \cup \W-= P_{\R}$, $\W+ \cap \W-= A_{\R} $,
and  define generalized Arnold surfaces
$\A{\pm}=\A{\pm}(f)=\W{\pm}\cup q(A)\subset Q$, where $Q=P/\!\conj$
and $q\:P\to Q$ is the quotient map.
If $A$ is connected and is of type  $1$ then
the Arnold surfaces $\A\pm_k\subset P$ can be defined 
like in the case $P=\Cp2$.

The collection of examples known to the author suggests
that the Rokhlin Conjecture might be  generalized as follows.

\tm{5.1. Generalized Rokhlin Conjecture}
If $(P,\conj)$ is a CDQ-surface
and $A\subset P$ is a Real, nonsingular, even
curve with $A_{\R}\ne\oo$
then the Arnold surfaces
in $P$ and in $Q$
are standard when  connected.
\endtm

By a {\it standard surface} in $Q$ we  mean a 
surface obtained as a
connected sum
of  standard surfaces in $S^4$ and possibly surfaces
$\pm(\Cp2, \text{conic})$,
$\pm(\Cp2, \varnothing)$,
$(S^2\times S^2, \varnothing)$.
We will consider multi-component standard surfaces
and admit connected sums in the both forms:
$(Y_1\#Y_2,F_1\#F_2)$ and $(Y_1\#Y_2,F_1\dsum F_2)$.
The above conjecture looks plausible for disconnected Arnold surfaces
as well, unless they include
some components of $X_{\R}$  entirely.
By a standard surface in $P$ we mean a connected sum of
$(P, C)$ and $(S^4, F)$,
where $C$ is a nonsingular complex curve in $P$
 and $F$  a standard surface in $S^4$.

Take $X=\{v\in L\, |\, v^2=f(p(v))\}$,
the double covering over $P$
with the projection $p|_{X}\: X\to P$,
and consider two conjugations
$\conj^{\pm}(v)=\conj_L(\pm v)$ which cover $\conj$.
Put $Y^{\pm}=X/\c{\pm}$.

\proclaim{5.2. Theorem }
If $(P,\conj)$ is a CDQ-surface and
the Arnold surface
$\A+$ $($or $\A-)$ for a Real curve $A$
is connected and standard, then
$(X,\conj^-)$ $($respectively $(X,\conj^+))$ is also CDQ-surface.
If $\A\pm$ is standard and has $k$ connected components, then
$Y^{\mp}\cong R\,\#\,(k\!-\!1)(S^1\!\times\! S^3)$,
where $R$ is a completely
decomposable 4-manifold.
\endproclaim

\demo{Proof}
If $\A{\pm}$ is connected then the proof is the same as for
Theorem 2.2 with only one additional remark
that the double branched
covering over $\pm(\Cp2, \text{conic})$ is a quadric
and, hence, standard.
If $\A{\pm}$ is not connected
we can use that
the double covering over
$(Y_1\#Y_2,F_1\dsum F_2)$ is obviously obtained from the
covering over
$(Y_1\#Y_2,F_1\#F_2)$
by adding  a 1-handle.
\qed
\enddemo
Theorem 5.3 below provides
 patterns of curves $A_i$ with standard Arnold surfaces.
Their deformations can produce further examples via Corollary 4.5.
Let $(P, \conj)$ be a nonsingular Real surface and
$(L, \conj_L)$ a real linear bundle of degree $d$. 
Consider a real pencil,
$f_t\: P\to L$,
 $f_t(x)=tf_1(x)+(1-t)f_0(x)$, $t\in \C$,
where $f_0, f_1$
are real holomorphic sections.
Denote by $D(f_t)$ the zero divisor of a section $f_t$ and assume 
that the curve $B=D({f_0})$ is nonsingular, connected
and intersects $D({f_1})$
transversally at real distinct points,
$b_1,\dots, b_d\in P_{\R}$.
Assume further that there is a real section $h\: P\to L$,
with $b_i\notin D(h)$, $1\le i\le d$.
Define sections
$v_{\e,t}, u_{\e,t}\: P\to L^{ 2}$,
parametrized by $t,\e\in\C$, as
\vskip-3mm
$$
v_{\e,t}=f_0^2+f_t^2-\e h^2,\ \
u_{\e,t}=-(f_0\!\cdot\! f_t+\e h^2).
$$
\vskip-1mm
It can be easily seen that the curves
$D(f_t)$, $D(v_{\e,t})$, $D(u_{\e,t})$ are
 Real, nonsingular and connected
for real and sufficiently small $\e,t>0$.
The real scheme of
$D(v_{\e,t})_{\R}$ of the curve $D(v_{\e,t})$
consists of $d$ small ovals around
the basepoints $b_i$.
 $D(u_{\e,t})_{\R}$ lies
between $D(f_0)$ and $D(f_t)$, inside $W^+(-f_0 f_t)$,
and is obtained from $B_{\R}$ by doubling
the  ovals of $B$ which does not contain basepoints
and replacing the ovals
 which contains $r$ basepoints by
$r$ ovals.

The following theorem asserts that these  curves have
standard Arnold surfaces  in $P$ and $Q$.

\proclaim{5.3. Theorem }
For  small enough $\e,t>0$
\roster
\item
$D({v_{\e,t}})$ is of type $1$
and has standard Arnold surfaces, $\A+_k(v_{\e,t})\subset P$,
isotopic to $B$; \ 
$D({u_{\e,t}})$ is of type  $1$ if $B$ is.
In the latter case
the Arnold surfaces $\A+_k(u_{\e,t})\subset P$
are standard and
isotopic to $B$.
\item
Arnold surfaces $\A+\subset Q$ for the both curves $u_{\e,t}$ and $v_{\e,t}$
bound  solid handlebodies embedded in $Q$.
These handlebodies are orientable in the case of
 $v_{\e,t}$ independently of the type of $B$.
In the case of
$u_{\e,t}$
they are orientable if and only if $B$ is of type 1.
\endroster
\endproclaim

\demo{Proof}
The degenerated curve $D({v_{0,t}})$
splits into a pair of conjugated imaginary curves,
$D(f_0\pm itf_1)$, isotopic to $B$.
After a perturbation they are deformed into a pair of halfs
of $D({v_{\e,t}})-D({v_{\e,t}})_\R$. 
$D({u_{0,t}})$ splits into a pair of Real curves, $B=D(f_0)$ and
$D(f_t)$, which are both of the same type for a small enough $t>0$.
If it is type 1, then the opposite halfs of $D(f_0)-D(f_0)_\R$
and $D(f_t)-D(f_t)_\R$ are fused after a perturbation
and form the halfs of
$D({u_{\e,t}})-D({u_{\e,t}})_\R$
(this way of gluing of ``halfs'' is well known 
and can be understood easily).
Consecutive degenerations as $\e\to 0$ and $t\to0$ provide isotopies
of the Arnold surfaces into $B$.
To prove (2) we use first the degeneration
 $\e\to 0$ and then
consider the trace of the variation $t\to 0$ as the solid handlebody.
In the case of $v_{\e,t}$
it gives an isotopy
between $\A+(v_{\e,t})$
and $\A+(v_{0,t})=q(D(f_0+itf_1))$.
The latter bounds the handlebody
$H'=q(H)$, where
 $H=\bigcup_{-t\leq\tau\leq t}D(f_0+i\tau f_1)$.
In the case of $u_{\e,t}$
 we get an isotopy between
$\A+(u_{\e,t})$ and $\A+(u_{0,t})$. The latter
bounds $H'=q(H)$, where
 $H=\bigcup_{0\leq\tau\leq t}D(f_\tau)$.

In the both cases $H'$ is smoothable at basepoints
since its links at $b_i$
are homeomorphic to discs
(see Figure 6).
Smoothability at the other points is obvious.
$H'$ is a solid handlebody, since it is homeomorphic
to the quotient of $A\times [-1,1]$ by the conjugation
$\tau_1(x,t)=(\conj|_A,\ t)$ in the case of $u_{\e,t}$,
and $\tau_2(x,t)=(\conj|_A,\ -t)$ in the case of $v_{\e,t}$.
\qed
\enddemo

Theorems 5.2 and 5.3 imply the following result.

\tm{5.4. Corollary}
Let $Q$ be simply connected and
 $Y^{\pm}$  defined as above via
$v_{\e,t}$ or  $u_{\e,t}$.
  Then
$Y^-\cong 2Q\#R$,
where $R$ is a completely decomposable $4$-manifold.

In particular, if $(P,\conj)$ is a CDQ-surface then
$(X,\conj^-)$ is CDQ as well.
\endtm

\rk{Remarks}
\roster
\item
If $\A+$ is orientable and $g\!=\!g(\A+)\!=\!g(B)$ is the genus,
then $R\cong g(S^2\times S^2)$.
The above condition is not satisfied if and only if
 $B$ is of type type 2 and we consider the case of $u_{\e,t}$.
Then  $R\cong g(\Cp2\,\#\,\barCP2)$.
\item
The proof of
Theorem 5.3 works also if $B$ is not connected
provided all components of $B$ have nonempty real part.
Then  $Y^-$ gets several additional $1$-handles,
the number of which is less by $1$ then the number of
connected components of $\A+$.
\item
In the case of $P=\Cp2$, the real part of curves $D(v_{\e,t})$
consists of $k^2$ disjoint ovals if $\deg(B)=k$.
Curves $D(u_{\e,t})$
have real schemes consisting of
$k^2$ disjoint ovals
added to
a doubled scheme
of a degree $k$ curve, which belongs to type 2 (and can be empty).
$L$-curves of degree $k$ can provide
real schemes of type 1 with $k^2$ empty ovals
as well.
It seems to be  interesting to study if
these three types of curves with $k^2$ ovals
are rigidly isotopic.
\endroster
\endrk

The arguments in the proof of Theorem 5.3 showing that
$\A+(v_{0,t})=q(D(f_0+itf_1))$ bounds a handlebody,
give the following

\tm{5.5. Corollary}
If $Q$ is simply connected,
 $z\in \C- \R$ and $|z|>0$ is small enough, then
$q(D(f_z))$
is a  standard surface in $Q$.
\endtm

\tm{5.6. Corollary}
Let  $C\subset P$ be a nonsingular
curve of degree $d$ which intersects
$P_{\R}$ at $d$ points. Then  $q(C)$ bounds a handlebody in $Q$
and, hence,
is a standard surface if $Q$ is simply connected.
\endtm

The latter Corollary follows from the previous one,
since the pencil of curves containing $C$ and $\conj(C)$ is real
and, therefore  contains Real curves.

The condition for the number of intersection points with
$P_{\R}$ is important; otherwise, $q(C)$ 
is not an embedded surface.

\tm{5.7. Corollary}
Nonsingular imaginary curves $A\subset \Cp2$
of degree $k$ which intersect  $\Rp2$ at $k^2$ points
are projected by the quotient map
into standard surfaces in $S^4=\Cp2/\!\conj$.
For example,  the images  of
imaginary  lines and
 imaginary conics intersecting $\Rp2$ in $4$ points
are unknots in $S^4$.
\endtm

\heading
\S 6 Fibered surfaces
\endheading

In this section we consider
complex surfaces
fibered over complex curves, $\pi\:P\to C$,
and study Arnold surfaces $\A{\pm}$ and quotients $Y^{\pm}$
associated with the
curves $A$ consisting of several nonsingular fibers.
A real structure on a fibered surface
is given, by definition, by a pair of commuting
complex conjugations
$\conj_C\: C\to C$
and $\conj\: P\to P$

$$
\CD
P     @>\conj>>     P\\
@V\pi VV              @V\pi VV\\
C      @>\conj_C>>       C\\
\endCD
$$

 Fibers over  real points $c\in C$, $F=\pi^{-1}(c)$, inherit real
structures, $\conj|_F$, from $P$.

We assume for the rest of the section that
\roster
\item
$P$, $C$ and a generic fiber $F$ are nonsingular and
connected.
\item
$P_{\R}\neq \oo$ and, in particular, $B_{\R}\neq \oo$.
\item
$A$ is  nonsingular, real and even.
Hence, it consists of
even number of real fibers
and several pairs of conjugated imaginary fibers
$ B_j=\pi^{-1}(b_j)$, $B_j'=\pi^{-1}(b_j')$, $b'_j=\conj(b_j)$,
$j=1,\dots, s$.
\item
Real fibers of $A$ are double, i.e., split into pairs
$A_i=\pi^{-1}(a_i),\ A_i'=\pi^{-1}(a_i')$,
where $ a_i, a_i'\in C_{\R}, i=1,\dots, r$,
are close enough to each other.
The latter means, in particular,
 that there are no singular real fibers between $a_i$ and $a_i'$.
\item
Real parts $A_{i,\R}$ of real fibers $A_i$, $i=1,\dots,r$ are not empty.
\endroster

Let $\W+_i\cong A_{i,\R}\times I$  denote the part of  $X_{\R}$ between
 $A_{i,\R}$ and $ A_{i,\R}',\ i=1,\dots, r$.
Let $q\:P\to Q$ be the quotient map and
$\A{\pm}$, $ X$, $\conj^{\pm}$ and $Y^{\pm}$  are defined
by the curve $A$ as in $\S 5$,
with the sign superscripts determined so that
$\A+=\bigcup_{i=1}^r\A+_i \bigcup_{j=1}^s{\goth B}_j$,
where  $\A+_i=\W+_i\cup q(A_i)\cup q(A_i')$ and
${\goth B}_j=q(B_j)=q(B_j')$.

\tm{6.1. Theorem}
$\A+$  bounds a disjoint union of $r+s$
solid handlebodies embedded in $Q$.
\endtm

Similar to Corollary 5.4 this implies the following

\tm{6.2. Corollary}
If $Q$ is simply connected then
\roster
\item
$\A+$ is standard,
\item
$Y^-\cong2Q\,\#\,R\,\#\,(r\!@!+\!@!s\!@!-\!@!1)(S^3\!\times\!
S^1)$, where $R$ is a completely
decomposable $4$-manifold.
\endroster
\endtm

\tm{6.3. Corollary}
If $(P,\conj)$ is a CDQ-surface and $Y^-$ is  simply connected,
then $(X,\conj^-)$ is CDQ.
\endtm

\demo{Proof of Theorem 6.1}
Let $[a_i,a_i']$ denote the small segment of $B_{\R}$ between $a_i$
and $a_i'$, and $H_i=\pi^{-1}([a_i,a_i'])$.
Then $H'_i=q(H_i)\cong q(A_i)\times[a_1,a_i']$
is a solid handlebody bounded by $\A+_i$.

Let $q_C\: C\to C/\!\conj_C$,
$\pi'\:Q\to C/\!\conj_C$  denote the natural projections.
Choose points $c_j\in(C_{\R}-\bigcup_{i=1}^r[a_i,a_i'])$, $j=1,\dots,s$,
with nonempty real parts of the fibers
$\pi^{-1}(c_j)$, for example  choose them near $a_1$ outside $[a_1,a_1']$.
Connect $q_C(b_j)$ with $c_j$ by disjoint set of
smoothly embedded arcs
$l_j\subset C/\!\conj_C$, which don't intersect with $C_{\R}$,
except  endpoints $c_j$,
and don't pass under singular fibers.
Then  $G'_j=(\pi')^{-1}(l_j)$ are
solid handlebodies bounded by ${\goth B}_j$.
$H_i'$ and $G'_j$ are disjoint by construction.
\qed
\enddemo

\rk{Remarks}
\roster
\item
Note that Theorem 6.1 implies Theorem 5.3(2), since
blow-ups at the basepoints of the pencil $f_t$
give a fibered surface and the curve $A$ is becoming
a double fiber. A blow-ups at a real basepoint does not change $Q$,
as it  is pointed out in the introduction.
A blow-up
 does not change also the Arnold surfaces $\A{\pm}_0$
of a curve $A_0$ when it is made at a cross-type point,
and it
 does not change  $\A+_0$ in the case of an isolated point of $A_0$.
Therefore  it does not change
 $Y^-$ in both cases as well.
\item
Orientability  of $H'$ and $G'$ is, clearly,
 determined by the same
criterion as in Theorem 5.3.
\endroster
\endrk

The topology of $\A-\subset Q$ is more complicated.
However, we can describe the double covering $Y^+\to Q$,
branched over $\A-$, in terms of some surgery on $P\to Q$.

\definition{Definition}
Given a pair of spaces $X_1, X_2$, with
 involutions $\theta_i\:X_i\to X_i,\ i=1,2$, we define
{\it $\Z/2$-product} of $(X_1,\theta_1)$ and $(X_2,\theta_2)$ as
the quotient space
\allowlinebreak
$X_1\times_{\Z/2}X_2=(X_1\times X_2)/\theta$,
by the product involution
 $\theta=\theta_1\times \theta_2$.
\enddefinition

Let $N_i=F\times_{\Z/2}\cil_i$, $i=1,2$, be $\Z/2$-products
of a real fiber $F=\pi^{-1}(a)$, $a\in C_{\R}$, supplied
with the involution $\conj|_F$,
and the cylinder $\cil_i=S^1\times[-1,1]$, supplied
with the involution $\tau_1(z,t)=(z,-t)$, for $i=1$,
and  with $\tau_2(z,t)=(-z,-t)$, for $i=2$.
Note, that
$\partial N_i=F\times_{\Z/2}(\partial \cil_i)\cong F\times S^1$.

\definition{Definition}
By {\it $\Z/2$-surgery} along a real nonsingular fiber
$F=\pi^{-1}(a)$
we mean replacing of
its tubular neighborhood
$\pi^{-1}(D)\cong D^2\times F$,
where $D\subset C$ is a small, regular, $\conj_C$-symmetric neighborhood
of $a$, by
$N_1$ ($\Z/2$-surgery of type 1),
or by $N_2$ ($\Z/2$-surgery of type 2).
\enddefinition

The diffeomorphism type of the result is well defined since
the trivializations of $\pi^{-1}(D)$ and $\partial N_i$
are, clearly, canonical up to diffeotopy.

The following examples are easy exercises.

\example\nofrills{Examples:}
\roster
\item
If  $F\cong S^2$, a rational fiber with $F_{\R}\neq \oo$,
then $\Z/2$-surgery of type 1 is the Morse modification of index 2
along $F$.
\item
If $F$ is a rational fiber and  $F_{\R}=\oo$,
then $\Z/2$-surgery of type 1
replaces $F\times D^2$ by $(\Rp3-D^3)\times S^1$.
\item
If $F$ is an elliptic fiber, $F_{\R}\ne \oo$, then
$\Z/2$-surgery of type 1 makes a logarithmic
transform of multiplicity $0$ along $F$,
which is, by definition, the modification
$D^2\times S^1\times S^1\to  S^1\times D^2\times S^1$,
that is
 the Morse  modification of index 1 in dimension 3,
multiplied by $S^1$.
\endroster
\endexample

\tm{6.4. Theorem}
 $Y^+$ can be obtained from $P$ by
$r$ $\Z/2$-surgeries of type $1$ along $A_i,\ 1\leq i\leq r$,
and $s$ $\Z/2$-surgeries of type $2$ along
arbitrary real fibers
$F_j\subset P,\ 1\le j\le s$, with $F_{j,\R}\ne \oo$.
\endtm

\demo{Proof}
Consider small regular $\conj$-symmetric neighborhoods
$\Cal H_i\supset H_i,\ \Cal G_j\supset G_j$ of the solid
handlebodies of Theorem 6.1.
The restriction of the covering $p\:X\to P$ over
$K=P-(\bigcup_i \Cal H_i\bigcup_j \Cal G_j)$ is trivial
and $\c+$ flips 2 copies of $K$. Thus  $p^{-1}(K)/\c+$
can be identified with $K$.

The restrictions of $\conj^+$ to the pull backs of neighborhoods,
$p^{-1}(\Cal H_i)$, $p^{-1}(\Cal G_j)$
are equivalent to $(\conj|_F)\times\tau_1$ on $F\times \cil_1$
and $(\conj|_F)\times\tau_2$ on $F\times \cil_2$.
Thus  quotients are $N_1, N_2$.
\qed
\enddemo

 Theorems 6.1---6.4 can be applied
 to elliptic surfaces; by
 making use of Example 3 we may summarise the above results as follows.

\tm{6.5. Theorem}
Let $P$ be an elliptic surface and $A=A_1\dsum A_1'$ a double real fiber.
Then
\roster
\item
$Y^+$ is diffeomorphic to the manifold,
obtained by
a logarithmic transform of multiplicity $0$ on $P$,
along the fiber $A_1$.
\item
If $Q$ is simply connected, then
$Y^-\cong 2Q\,\#\, R$, where
$R$ is
$ S^2\times S^2$ if $A_{1,\R}$ belongs to type  1,
and $\Cp2\#\,\barCP2$ if it belongs to type  2.
\newline
In particular, if $(P,\conj)$ is CDQ-surface then  $(X,\conj^-)$ 
is CDQ as well.
\endroster
\endtm

\tm{6.6. Corollary}
Let $P=E(n)_{m_1,\dots, m_k}$. Then $Y^+$ is completely decomposable if
it is simply connected.
\endtm

The last corollary follows from Theorem 6.5,
since $Y^+\!=E(n)_{m_1,\dots, m_k,0}$ is not simply connected
unless $k=0$; however $E(n)_0$ is decomposable (see \cite{G}).

\heading
\S 7. Some applications
\endheading

In this section we consider
applications of the results of \S 3--5
to Real curves $A\subset\Cp2$ of degree  $\le 6$
and to Real $K3$ surfaces.

All nonsingular curves $A$ of degree $\le 4$ with $A_{\R}\ne \oo$
are known to be $L$-curves \cite{F1}.
Hence, Theorem 3.1 implies

\tm{7.1. Corollary}
Arnold surfaces $\A{\pm}(A)$ for nonsingular
curves $A$ of degree $2$ and $4$
with $A_{\R}\ne\oo$
are standard.
\endtm

The case $\deg(A)=6$ is a somewhat more delicate.
According to the classification of sextics (see, e.g., \cite{V})
$A$ is determined up to rigid isotopy by the arrangement
of the components ({\it ovals}) of $A_{\R}$ in $\Rp2$
and the type of $A$.
Similar to \cite{V},
we will use the notation
$\shhh\a\b \varkappa$ to code  curves of type $\varkappa=1,2$
which have
an oval containing $\b$ separate ovals inside and 
$\a$ separate ovals outside.
In the case of sextics the total number of ovals is $\a+\b+1\le 11$.
In addition to the schemes $\shhh\a\b \varkappa$
sextics may have schemes $\ss\a \varkappa$, $\a\le 10$,
corresponding to $\a$ separate ovals,
and the scheme $\sss$, a nest of 3 ovals.
We may omit the type
$\varkappa$ in codes of the schemes
when it  is determined
by the arrangement of ovals of a sextic,
like we did it for the above nest, which can be only of type 1
(for possible schemes and their types see \cite{V}).

We will also use notation $\shhhh\a\b \varkappa{\pm}$ and likewise 
for the domain $\W{\pm}(A)$ of a curve $A$ 
with the real scheme $\shhh\a\b \varkappa$.

\tm{7.2. Theorem}
All nonsingular Real sextics $A$ with $A_{\R}\ne\oo$
have standard Arnold surfaces $\A+(A)$.
If the scheme of $A_{\R}$ differs from
$\shhh 191$, $\shhh182$,
 $\la1\la9\ra\ra_2$ and
 $\la1\la8\ra\ra_2$, then
it has also a standard
Arnold surface  $\A-(A)$.
\endtm

\demo{Proof}
It is well known and
not difficult to see directly from the
Hilbert and Gudkov constructions
of nonsingular Real sextics (cf.~\cite{V}),
that the ones with schemes
$\shhh\a\b \varkappa,\, \varkappa=1,2$, can be
deformed to the both schemes
\roster
\item
 $\shhh{(\a-1)}\b 2$ if $\a\ge 1$, and
\item
 $\shhh\a{\b-1}2$ if $\b\ge2$, or to
$\ss\a2$ if $\b=1$
\endroster
by passing through a cross-like real node
which connects the ambient oval with an exterior oval
in the first deformation and with an interior one in the second.
The only exception is
the scheme $\shh91$, 
which can be reduced 
to $\ss{10}2$ 
only by 
contracting the interior oval.

Figure 8 shows that sextics
 with the nest scheme $\sss$ and with the scheme $\ss{10}2$ are
$L$-curves, thus  have standard Arnold surfaces.
\midinsert
\topcaption{Figure 8}
Construction of $L$-curves $\sss$ and $\ss{10}2$
\endcaption
\vspace{50mm}
\endinsert

The latter scheme gives the maximal possible value
of $\chi(W^{\pm})$ realizable for sextics, $\chi(W^+)=10$.
We will show that
the other schemes $W^\pm$ can be
 obtained  from it by moves $M_1$,
with 4 exceptions mentioned in Theorem 7.2.
Hence, Corollary 4.5 can be applied
to provide the result of the theorem.

It follows from the above remark on reduction of schemes that
$\ss{10}2^+\red\ss{\a}2^+$ and
$\ss{\a+1}2^+\red\shhhh\a\b \varkappa+$
for all schemes with
$\a\,\ge 0,\b\ge1,\ \varkappa=1,2$,
whenever they are realizable by sextics,
except  $\shh91^+$. 
Further, by Theorem 5.3 sextics
with schemes $\ss91$, $\shh91$
have standard $\A+$ as well.
Altogether this shows that $\A+$ are standard in all the cases.
Figure 9 shows that $\ss92^+\red\la1\la8\ra\ra^-_1$.
By the same arguments as above we get
 $\la1\la8\ra\ra^-_1
\red\la1\la\b\ra\ra^-_2\red
\shhhh\a\b \varkappa-$
for all existing schemes of sextics with
 $ 0\le\a,\ 0\le\b\le7, \ \varkappa=1,2$.

This includes all cases with the $4$ exceptions mentioned in
Theorem 7.2.
\qed
\enddemo
\midinsert
\topcaption{Figure 9}
 Deformation
$\ss92^+\red\la1\la8\ra\ra^-_1$.
The middle curve on the top is obtained by the Hilbert method via
the perturbation  shown on the bottom
(for details on Hilbert curves
see \cite{V})
\endcaption
\vspace{65mm}
\endinsert

\tm{7.3. Theorem}
Any Real nonsingular $K3$-surface $(X,\conj)$ with $X_{\R}\ne\oo$
is CDQ-surface.
\endtm

\demo{Proof}
The components of the moduli space for nonsingular Real $K3$
surfaces $(X,\conj)$ 
can be described in terms of arithmetics in $H_2(X)$ and
distinguished  by 
topological types of $X_{\R}$ and
vanishing of $[X_{\R}]\in H_2(X;\Z/2)$.
This fact is well known and
follows from the global Torelli theorem, 
as it was noticed by V. Kharlamov and V. Nikulin \cite{N}.
(In fact, \cite{N} considers
only algebraic K3 surfaces, but after Torelli theorem was
 extended to abstract K3 
the arguments of \cite{N} concerns the latter case as well.)
It is known also that
$X_{\R}$ may be homeomorphic to $S_1\dsum S_1$, or $S_g\dsum kS_0$,
$g+k\le 11$, with certain restrictions to $g,k$
when $g+k=11$ or $10$
(cf.~\cite{N} or \cite{V});
here  $S_g$ denotes a genus $g$ orientable surface 
and $kS_0$ a disjoint union of $k$ spheres.

After comparing the classification of Real $K3$-surfaces
with the classification of Real sextics, cf.~\cite{N, V},
we see
that all but one deformation types of
 Real $K3$-surfaces  are represented by double coverings over $\Cp2$
branched along sextics, as $(X, \conj^{\pm})$ in the notation of \S 2.
More precisely,
 all the types of
$K3$ surfaces for which $X_{\R}$ has a non sphere component
can be represented as $(X, \conj^-)$.
The ones for which  $X_{\R}\cong kS_0$
and $[X_{\R}]\ne 0\in H_2(X;\Z/2)$
can be represented as $(X, \conj^+)$,
when $A$ is
a sextic having real scheme $\ss k{}$.

By Theorem 7.2,
  Arnold surfaces $\A+$ are standard for sextics with any scheme,
and $\A-$ for sextics with schemes $\ss k{}$.
Hence, by Theorem 2.2
$(X,\conj)$ are CDQ-surfaces
in all the cases
(recall that $\A{\pm}$ corresponds to $\conj^{\mp}$).

The only real deformation type of $K3$ surfaces, $(X,\conj)$,
which does not contain
a double plane
has $X_{\R}\cong 8S_0$ and
$[X_{\R}]=0$.
We will see that there is a deformation of  $(X,\conj)$
fusing a pair of components of $X_{\R}$
and giving as the result  $(X',\conj')$ with
 $X'_{\R}=7S_0$.
$(X',\conj')$ is already shown to be CDQ-surface, hence, by
Corollary 4.4 $(X,\conj)$ is CDQ-surface as well.
To see the fusing deformation connecting
 $(X,\conj)$ with $(X',\conj')$
consider a family of
double quadrics branched along Real curves of degree $(4,4)$.
It consists of the Real $K3$ surfaces and includes
$(X,\conj)$ not realizable as a double plane.
Explicitly,
take the quadric $P\subset\Cp3$ defined by the equation
$x^2+y^2-z^2=t^2$, the projection $f\: P\to \Cp2$,
$(x:y:z:t)\mapsto (x:y:z)$,
and the pull back $\widetilde A=f^{-1}(A)\subset P$ of a
plane curve $A$ of degree 4, which real part $A_{\R}$
has 4 ovals outside the oval of
the branched locus of $f$, i.e.
the conic 
$x^2+y^2-z^2=0$.
$\widetilde A_{\R}$ is of type 1 and
has 8 separate ovals lying in the hyperboloid $P_{\R}$.
$(X,\conj)$ is
 the double covering over $P$, branched along $\widetilde A$.
A deformation
 fusing a pair of components of $X_{\R}$ is obtained then
from a deformation of  $\widetilde A$, which fuses a pair of
ovals.
\qed
\enddemo

Standard calculations (see Supplement below) give values
$b^+_2(Y)=1, \ \bn-=9+\frac12\chi(X_{\R})$ for
Real $K3$ surfaces $(X,\conj)$.
Hence, by  Theorem 7.3 either
$Y\cong\Cp2\,\#\,k\barCP2$, where $k=9+\frac12\chi(X_{\R})$,
or  $Y=S^2\times S^2$.
The latter holds
in 2 cases:
\roster
\item
$X_{\R}\cong S_{10}\dsum S_0$,
\item
$X_{\R}=S_9$ and the class
$[X_{\R}]\in H_2(X;\Z/2)$ vanishes,
\endroster
as it follows from Theorem 5.3, since
these two types of $K3$ surfaces are represented by double planes
$(X,\conj^-)$ branched along sextics with one of the schemes
$\shh91$,  $\ss91$,
and because
$Y^\mp=S^2\times S^2$ if and only if
the Arnold surface $\A\pm$ is orientable.

\rk{{\bf 7.4.} Supplement}
We reproduce here the calculation
of $\bn{\pm}$ in more general setting, which
completes the description of $Y$ if it is decomposable.
By the Hirzebruch signature formula
and the Riemann--Hurwitz formula, we have
$$
\gathered
\sigma(X)=
2\sigma(Y)-X_{\R}\circ X_{\R},\\
\chi(X)=2\chi(Y)-\chi(X_{\R}),
\endgathered
$$
where $\sigma$ denotes signature and the self-intersection number
$X_{\R}\circ X_{\R}$
is taken in $X$. We have also
$X_{\R}\circ X_{\R}=-\chi(X_{\R})$,
as it is pointed out in the proof of Corollary 4.8.
This gives if $b_1(X)=0$ and, hence, $b_1(Y)=0$
$$
\gathered
\bn+=\tfrac12(\bp{X}-1),\\
\bn-=\tfrac12(\bm{X})+\chi(X_\R)-1),
\endgathered
$$
One can also eliminate terms involving $X$
when
 $(X,\conj^\pm)$ is obtained as the double covering
 over $(P,\conj)$,
 branched along curve $A$,
by making use of the relations
$$
\gathered
\sigma(X)=2\sigma(P)-A\circ A,\\
\chi(X)=2\chi(P)-\chi(A).
\endgathered
$$
We can, express $b_2^\pm(Y^\pm)$
in terms of the zero  divisor $B$
of $L$ (for notation see \S 5), by making use of adjunction formula,
 applied  to $A$ and $B$ in $P$.
$$
\gathered
\bp {Y^\pm}=\bp P-\tfrac12\chi(B)=\bp P+\tfrac12(B+K)\!\cdot\! B,\\
\bm {Y^\pm}=\bm P+\tfrac12(-\chi(B)+\chi(X^\pm_{\R}))+d=
\bm P +\tfrac12(3B+K)\!\cdot\! B+\tfrac12\chi(X^\pm_{\R}),
\endgathered
$$
where $d=B\!\cdot\! B=\deg(L)$ and $K$ is the canonical divisor on $P$.
\endrk

\enddocument